 \documentclass[preprint,review,12pt]{elsarticle}

\usepackage{amssymb}

\usepackage{amsmath}
\usepackage{graphicx}
\usepackage{graphics}
\usepackage{subfigure}
\usepackage{mathptmx}
\usepackage[misc]{ifsym}
\usepackage{algorithm}
\usepackage{algorithmicx}
\usepackage{algpseudocode}
\usepackage{graphicx}
\usepackage{graphics}
\usepackage{subfigure}
\usepackage{epsfig}
\usepackage{color}
\usepackage{multirow}

\usepackage[top=2cm, bottom=2cm, left=3cm, right=3cm] {geometry}

\newtheorem{theorem}{Theorem}
\newtheorem{lemma}{Lemma}

\newproof{proof}{Proof}

\newtheorem{assumption}{Assumption}

\begin{document}

\begin{frontmatter}



\title{{\bf Distributed Convex Optimization With Coupling Constraints Over Time-Varying Directed Graphs}
}


\author[cq]{Chuanye Gu} 
\author[cq]{Zhiyou Wu}  
\author[cq]{Jueyou Li \corref{cor1}} \ead{lijueyou@163.com}
\author[cq]{Yaning Guo} 
 \cortext[cor1]{}

\address[cq]{School of Mathematical Sciences, Chongqing Normal University, Chongqing, 400047, China}

\begin{abstract}
This paper considers a distributed convex optimization problem over a time-varying multi-agent network, where each agent has its own decision variables that should be set so as to minimize its individual objective subject to local constraints and global coupling equality constraints. Over directed graphs,  a distributed algorithm is proposed that incorporates the push-sum protocol into dual subgradient methods. Under the convexity assumption, the optimality of primal and dual variables, and constraint violations is first established. Then the explicit convergence rates of the proposed algorithm are obtained. Finally, some numerical experiments on the economic dispatch problem are provided to demonstrate the efficacy of the proposed algorithm.
\end{abstract}

\end{frontmatter}


\section{Introduction}\label{sec1}

Due to the emergence of large-scale networks, distributed optimization problems have attracted recently considerable interests in many fields such as the control and operational research communities, wireless and social networks \cite{baingana2014proximalgradient}, \cite{mateos2012distributed}, power systems \cite{bolognani2015distributed}, \cite{zhang2016distributed}, robotics \cite{martinea2007on}, and name a few. These problems share some common characteristics: the entire optimization objective function can be decomposed into the sum of several individual objective functions over a network, and each individual only knows its own objective function, and only individual and its neighbors cooperate to solve the problem by interacting the network information locally \cite{nedic2009distributed}. Many researchers have investigated various multi-agent optimization problems arising in the engineering community \cite{gharesifard2012distributed,lorenzo2016next,jakovetic2014fast,li2015gradient}.

In literatures, consensus based distributed algorithms for solving distributed problems are mainly divided into three classes: primal consensus based algorithms, dual consensus based algorithms and primal-dual consensus based algorithms, see \cite{nedic2009distributed,duchi2012dual,zhu2012on,Shi2014On,Shi2015EXTRA}
In \cite{nedic2009distributed}, Nedi\'{c} et al. firstly proposed a distributed subgradient algorithm  with provable convergence rates, while its stochastic variant was investigated in \cite{ram2010distributed} and its asynchronous variant in \cite{nedic2011asynchronous}. Based on dual averaging methods, Duchi et al. \cite{duchi2012dual} proposed a distributed dual averaging algorithm and obtained the convergence rate, scaling inversely in the spectral gap of the networks. Resorting to Lagrange dual method, the papers \cite{zhu2012on,yuan2012distribted} designed a distributed primal-dual algorithm for solving distributed problem with equality or inequality constraints. However, distributed methods proposed in most to previous works require the use of doubly stochastic weight matrices, which are not easily constructed in a distributed fashion when the graphs are directed.

To overcome this issue, the work in \cite{tsianos2012push} proposed a different distributed subgradient approach in directed and fixed network topology, in which the messages among agents is propagated by ``push-sum" protocol. The push-sum based distributed method eliminates the requirement of graph balancing, however, the communication protocol is required to know the number of agents or the graph. Although the authors in \cite{nedic2015distributed} canceled the requirement of a balanced graph and proposed a push-subgradient approach with an explicit convergence rate at order of $\mathcal{O}(\ln t/\sqrt{t})$, they only investigated the unconstrained  distributed optimization  problems. Very recently, the reference \cite{nedic2017Achieving} proposed a Push-DIGing method that uses column stochastic matrices and fixed step-sizes, which can achieve a geometric rate.

The problems for solving distributed optimization subject to equality or (and) inequality constraints has received considerable attentions \cite{zhu2012on,yuan2012distribted,chang2015Multi-agent,maegellos2016distributed,Yuan2016On,Yuan2016Regularized}. In \cite{zhu2012on}, Zhu et al. firstly proposed a distributed Lagrangian primal-dual subgradient method by characterizing the primal-dual optimal solutions as the saddle points of the Lagrangian function related to the problem under consideration. Yuan et al. in \cite{yuan2012distribted} developed a variant of the distributed primal-dual subgradient method by introducing multistep consensus mechanism. For more general distributed optimization problems with inequality constraints that couple all the agents' decision variables, Chang et al. \cite{chang2015Multi-agent} designed a novel distributed primal-dual perturbed subgradient method and obtained the estimates on convergence rate, also see \cite{chang2016proximal}. By making use of dual decomposition and proximal minimization,  Falsone et al. \cite{falsone2016dual} proposed a novel distributed method to solve inequality-coupled optimization problems. Their proposed approach
can converge to some optimal dual solution of the centralized problem counterpart, while the primal variables converge to the set of optimal primal solutions.
However, the implementation of the algorithm proposed in \cite{falsone2016dual} requires the double stochasticity of communication weight matrices, without any estimates on the convergence rate of their proposed method.

In this paper, we investigate a distributed optimization problem subject to coupling equality constraints over time-varying directed networks. Under the framework of dual decomposition, we propose a distributed dual subgradient method with push-sum protocol to solve this problem. Under the assumption of directed graphs, we prove the optimality of dual and primal variables, and obtain the explicit convergence rates of the proposed method.

Compared to existing literatures, the contributions of this paper are two folds:

i) Relaxations undirected graphs to directed graphs. The work in \cite{nedic2015distributed} proposed a push-sum based distributed method to unconstrained optimization problems. By resorting to dual methods and push-sum protocols, we propose distributed dual subgradient method to solve a class of convex optimization subject to coupling equality constraints. Our algorithm can be viewed as an extension of push-sum based algorithms \cite{nedic2015distributed} to a constrained setting.
The consensus based primal-dual distributed methods proposed in \cite{zhu2012on,yuan2012distribted} require that the networks are undirected and the communication weight matrices are double stochastic, which are unrealistic over directed networks. By utilizing the push-sum scheme considered in \cite{tsianos2012push,nedic2015distributed}, our method can deal with distributed optimization problems over time-varying directed graphs, only needing the column stochastic matrices.

ii) Estimates on the convergence rate of the proposed method. The reference in \cite{falsone2016dual} analyze the convergence of dual optimality and primal optimality, without investigating constraint violations of the problem interest. In our algorithm, we extend the algorithm in \cite{falsone2016dual} to directed graphs, and establish the convergence results for dual optimality, primal optimality and constraint violations. More importantly, we obtain the explicit convergence rate of the proposed algorithm.

The remainder of this paper is organized as follows. We state the problem and related assumptions in Section \ref{sec2}. Section 3 proposes the solution method and main results. Section 4 provides the proof of main results. Numerical simulations are given in Section 5. Finally, Section 6 draws some conclusions.

Notation: We use boldface to distinguish the scalars and vectors in $\mathbb{R}^{n}$. For instance, $\nu_{i}[t]$ is a scalar and $\mathbf{u}_{i}[t]$ is a vector. For a matrix $D$, we will use the $D_{ij}$ to show its $i,j$'th entry. We use $||\mathbf{x}||$ to denote the Euclidean norm of a vector $\mathbf{x}$, $||\mathbf{x}||_{1}$ denote the $\ell_{1}$ norm of a vector $\mathbf{x}$, and $\mathbf{1}$ represent the vector of ones.

\section{ Distributed optimization with coupling equality constraints} \label{sec2}

\subsection{Problem statement and dual decomposition}

Consider a time-varying network with $m$ agents which would like to cooperatively solve the following minimization problem:
\begin{equation}\label{cop}
   \min_{\{\mathbf{x}_{i}\in \mathbf{X}_{i}\}_{i=1}^{m}} ~~F(\mathbf{x}):=\sum_{i=1}^{m} f_{i}(\mathbf{x}_{i}) ~~~~\textrm{s.t.}~~~~ \sum_{i=1}^{m}(A_{i}\mathbf{x}_{i}-\mathbf{b}_{i})=0,
\end{equation}
where each agent $i,i=1,\ldots,m$ only knows its own vector $\mathbf{x}_{i}\in \mathbb{R}^{n_{i}}$
of $n_{i}$ decision variables, its local constraint set $\mathbf{X}_{i}\subseteq \mathbb{R}^{n_{i}}$, objective function $f_{i}(\mathbf{x}_{i}): \mathbb{R}^{n_{i}}\rightarrow \mathbb{R}$, and all agents subject to the coupling equality constraints $\sum_{i=1}^{m}(A_{i}\mathbf{x}_{i}-\mathbf{b}_{i})=0$, $A_{i}\in \mathbb{R}^{p\times n_i}$ and $\mathbf{b}_{i}\in \mathbb{R}^{p}$. $\mathbf{x}=(\mathbf{x}_{1}^{\top},\mathbf{x}_{2}^{\top},\cdots, \mathbf{x}_{m}^{\top})^{\top}$ with $n=\sum_{i=1}^{m} n_{i}$, belongs to $\mathbf{X}=\mathbf{X}_{1}\times \mathbf{X}_{2}\times\cdots\times \mathbf{X}_{m}$.

Problem (\ref{cop}) is quite general arising in diverse applications. For examples, distributed model predictive control\cite{Wang2017distributed}, network utility maximization\cite{beck2014}, real-time pricing problems for smart grid \cite{zhang2012convergence,zhang2016distributed,li2016a} can be modeled in this class of problems.

To decouple the coupling equality constraints, we utilize the Lagrange dual method.
Firstly, we introduce the Lagrangian function $L(\mathbf{x},\lambda)$ of problem (\ref{cop}), given by
\begin{equation}\label{lf}
   L(\mathbf{x},\lambda) = \sum_{i=1}^{m}\{f_{i}(\mathbf{x}_{i})+ \lambda^{\top}(A_{i}\mathbf{x}_{i} - \mathbf{b}_{i})\}= \sum_{i=1}^{m}L_{i}(\mathbf{x}_{i},\lambda),
\end{equation}
where $L(\mathbf{x},\lambda)$: $\mathbb{R}^{n}\times \mathbb{R}^{p}\rightarrow \mathbb{R}$, $\lambda\in \mathbb{R}^{p}$ is the vector of Lagrange multipliers.
Define the dual function of problem (\ref{cop}) as follows
$$\phi(\lambda):= \min_{\mathbf{x}\in \mathbf{X}}L(\mathbf{x},\lambda).$$
Note that the Lagrangian function $L(\mathbf{x},\lambda)$ is separable with respect to $\mathbf{x}_{i}, i=1,\ldots,m$. Thus, the dual function $\phi(\lambda)$ can be rewritten as
\begin{equation}\label{li}
  \phi(\lambda)= \sum_{i=1}^{m} \phi_{i}(\lambda)
  =\sum_{i=1}^{m}\min_{\mathbf{x}_{i}\in \mathbf{X}_{i}}L_{i}(\mathbf{x}_{i},\lambda),
\end{equation}
where $\phi_{i}(\lambda)$ can be regarded as the dual function of agent $i, i=1,\ldots,m$. It is obvious that the dual function $\phi(\lambda)$ is concave but non-smooth generally.

Then, the dual problem of problem (\ref{cop}) can be written as $\max_{\lambda}\min_{\mathbf{x}\in X} L(\mathbf{x},\lambda)$, or, equivalently,
\begin{equation}\label{dual}
 \max_{\lambda}\sum_{i=1}^{m} \phi_{i}(\lambda).
\end{equation}
The coupling equality constraints between agents is represented by the fact that $\lambda$ is a common decision vector and all the agents should agree on its value.

\subsection{Assumptions }

The following assumptions on the problem (\ref{cop}) and on the communication time-varying network are required to show properties of convergence for the proposed algorithm.

\begin{assumption}\label{A1}
For each $i=1,2,\ldots, m$, the function $f_{i}(\cdot)$: $\mathbb{R}^{n_{i}}\rightarrow \mathbb{R}$ is convex, and the set $\mathbf{X}_{i}\subseteq \mathbb{R}^{n_{i}}$ is non-empty, convex and compact.
\end{assumption}

Note that, under Assumption \ref{A1}, for any $\mathbf{x}_{i}\in \mathbf{X}_{i}$, there is a constant $G_{i}>0$ such that $||A_{i}\mathbf{x}_{i}-\mathbf{b}_{i}||\leq G_{i}$, due to the compactness of $\mathbf{X}_{i}$, $i=1,2,\ldots, m$. Let $G=\sum_{i=1}^{m}G_{i}$.

\begin{assumption}\label{A2}
The Slater's condition of problem (\ref{cop}) holds, i.e., there exists $\mathbf{\check{x}}=(\mathbf{\check{x}}_{1}^{\top},\ldots, \mathbf{\check{x}}_{m}^{\top})^{\top}\in relint(X)$ such that $\sum_{i=1}^{m} (A_{i}\mathbf{\check{x}}_{i}-b_{i})=0$, where $relint(X)$ is the relative interior of the constrained set $X$.
\end{assumption}

Under Assumptions \ref{A1} and \ref{A2}, the strong duality holds and an optimal primal-dual pair $(x^{*},\lambda^{*})$ exists \cite{boyd2004convex}, where $\mathbf{x}^{*}=(\mathbf{x}_{1}^{*\top},\ldots, \mathbf{x}_{m}^{*\top})^{\top}\in \mathbb{R}^{n}$ and $\mathbf{\lambda}^{*}\in \mathbb{R}^{p}$ are optimal solutions of the primal problem (\ref{cop}) and dual problem (\ref{dual}), respectively. Moreover, the saddle-point theorem also holds \cite{boyd2004convex}, i.e., given an optimal primal-dual pair $(\mathbf{x}^{*},\mathbf{\lambda}^{*})$, we have that
\begin{eqnarray}\label{ap}
L(\mathbf{x}^{*},\lambda)\leq L(\mathbf{x}^{*},\lambda^{*})\leq L(\mathbf{x},\lambda^{*}),~\forall \mathbf{x}\in \mathbf{X},\lambda\in \mathbb{R}^{p}.
\end{eqnarray}
Let $X^{*}$ and $\Lambda^*$ be the optimal solution set of the primal problem (\ref{cop}) and dual problem (\ref{dual}), respectively.

We assume that each agent can communicate with other agents over a time-varying network. The communication topology is modeled by a \emph{directed} graph $\mathcal{G}[t]=(\mathcal{V}, \mathcal{E}[t])$ over the vertex set $\mathcal{V}=\{1,\ldots,m\}$ with the edge set $\mathcal{E}[t]\subseteq \mathcal{V}\times \mathcal{V}$.
Let $\mathcal{N}_{i}^{in}[t]$ represent the collection of in-neighbors and $\mathcal{N}_{i}^{out}[t]$ represent the collection of out-neighbors of agent $i$ at time $t$, respectively. That is,
$$\mathcal{N}_{i}^{in}[t]:=\{j|(j,i)\in \mathcal{E}[t]\}\cup  \{i\},$$
$$\mathcal{N}_{i}^{out}[t]:=\{j|(i,j)\in \mathcal{E}[t]\}\cup  \{i\},$$
where $(j,i)$ represents agent $j$ may send its information to agent $i$.
And let $d_{i}(t)$ be the out-degree of agent $i$, i.e.,
$$d_{i}[t]=|\mathcal{N}_{i}^{out}[t]|,$$

We introduce a time-varying communication weight matrix $D[t]$ with elements $(D[t])_{ij}$, defined by
\begin{eqnarray}\label{cop2}
& (D[t])_{ij} =  ~~~~\left\{
\begin{aligned}
&~~~~\frac{1}{d_{j}[t]},~~\textrm{when}~ j\in \mathcal{N}_{i}^{in}[t],~i,~j=1,2,\ldots,m,\\
&~~~~~~0,~~~~~~~~\textrm{otherwise}.\\
\end{aligned}
\right. 
\end{eqnarray}

Note that the communicated weight matrix $D[t]$ is column-stochastic. In our paper, we do not require the assumption of double-stochastic on $D[t]$.
We need the following assumption on the weight matrix $D[t]$, which can be found in \cite{nedic2015distributed}, \cite{nedic2010constaints}.
\begin{assumption}\label{A3}
i) Every agent $i$ knows its out-degree $d_{i}[t]$ at every time $t$;~
ii) The graph sequence $\mathcal{G}[t]$ is $B$-strongly connected, namely, there exists an integer $B>0$ such that the sequence $\mathcal{G}[t]$ with edge set $\mathcal{E}[t]=\cup_{l=kB}^{(k+1)B-1}\mathcal{E}[l]$ is strongly connected, for all $k\geq0$.
\end{assumption}

\section{Algorithm and main results}\label{sec3}

\subsection{Distributed dual sub-gradient push-sum algorithm}

Generally, the problem (\ref{cop}) could be solved in a centralized manner \cite{boyd2004convex}. However, if the number $m$ of agents is significantly large, this may cause computational challenge. Additionally, each agent would
be required to share its own information, such as the objective $f_{i}$, the constraints $X_{i}$ and $(A_{i},\mathbf{b}_{i})$, either with the other agents or with a central coordinate collecting
all information, which is possibly undesirable in many cases, due to privacy concerns \cite{falsone2016dual}.

To overcome both the computational challenge and the privacy issues stated above, we propose a  Distributed Dual Sub-Gradient Push-Sum algorithm (DDSG-PS, for short) by resorting to solve the dual problem (\ref{dual}). Our proposed algorithm DDSG-PS is motivated by the sub-gradient push-sum method \cite{nedic2015distributed} and dual decomposition \cite{li2016a,falsone2016dual}, described as in Algorithm \ref{alg1}.

\begin{algorithm}[!htb]
\caption{Distributed Dual Sub-Gradient Push-Sum Algorithm (DDSG-PS) }\label{alg1}
\begin{algorithmic}[1]
\State Initialization: for $i=1,2,\ldots,m$, given $\mathbf{\mu}_{i}[0]\in \mathbb{R}^{p}$ and $\nu_{i}[0]=1$; set $t:=0$;
\Repeat
\For{each agent $i=1,\ldots,m$}
\State $\mathbf{u}_{i}[t+1]=\sum_{j=1}^{m}(D[t])_{ij}\mathbf{\mu}_{j}[t]$;
\State $\nu_{i}[t+1]=\sum_{j=1}^{m}(D[t])_{ij}\nu_{j}[t]$;
\State $\mathbf{\lambda}_{i}[t+1]=\frac{\mathbf{u}_{i}[t+1]}{\nu_{i}[t+1]}$;
\State $\mathbf{x}_{i}[t+1]=\arg\min_{\mathbf{x}_{i}\in \mathbf{X}_{i}}\{f_{i}(\mathbf{x}_{i})
+\mathbf{\lambda}_{i}[t+1]^{\top}(A_{i}\mathbf{x}_{i}-\mathbf{b}_{i})\};$
\State $\mathbf{\mu}_{i}[t+1]=\mathbf{u}_{i}[t+1]+\beta[t+1](A_{i}\mathbf{x}_{i}[t+1]-\mathbf{b}_{i})$;
\EndFor
\State set $t=t+1$;
\Until{a preset stopping criterion is met} 
\end{algorithmic}
\end{algorithm}

In Algorithm \ref{alg1}, each agent $i$ broadcasts (or pushes) the quantities $\mathbf{\mu}_{i}[t]/d_{i}[t]$ and $\nu_{i}[t]/d_{i}[t]$ to all of the agents in its out-neighborhood $\mathcal{N}_{i}^{out}[t]$. Then, each agent simply sums all the received messages to obtain $\mathbf{u}_{i}[t+1]$ in step 4 and $\mathbf{u}_{i}[t+1]$ in step 5, respectively. The update rules in steps 6-8 are implemented locally. In particular, the update of local primal vector $\mathbf{x}_{i}[t+1]$ in step 7 is performed by minimizing $L_{i}$ with respect to $\mathbf{x}_{i}$ evaluated at $\lambda=\lambda_{i}[t+1]$, while the update of the dual
vector $\mathbf{\mu}_{i}[t+1]$ in step 8  involves  the
maximization of $L_{i}$ with respect to $\lambda_{i}$ evaluated at $\mathbf{x}_{i}=\mathbf{x}_{i}[t+1]$.
Note that the term $A_{i}\mathbf{x}_{i}[t+1]-\mathbf{b}_{i}$ in step 8 is the sub-gradient of the dual function $\phi_{i}(\lambda)$ at $\lambda=\lambda_{i}[t+1]$.

It is shown in \cite{nedic2010constaints,beck2014} that the local primal vector $\mathbf{x}_{i}[t]$ does not converge to the optimal solution $\mathbf{x}_{i}^{*}$ to problem (\ref{cop}) in general. As compared to $\mathbf{x}_{i}[t]$, however, the following recursive auxiliary primal iterates
$$\widehat{\mathbf{x}}_{i}[t+1]=\widehat{\mathbf{x}}_{i}[t]+\frac{\beta[t+1]}{\sum_{r=1}^{t+1}
\beta[r]}(\mathbf{x}_{i}[t+1]-\widehat{\mathbf{x}}_{i}[t]), ~\mathrm{for~all}~t\geq 0$$
shows better convergence properties with $\widehat{\mathbf{x}}_{i}[0]=\mathbf{x}_{i}[0]$ \cite{nedic2009distributed,zhu2012on,chang2015Multi-agent}. In a similar way, we introduce an recursive auxiliary dual iterates as follows
$$\widehat{\mathbf{\lambda}}_{i}[t+1]=\widehat{\mathbf{\lambda}}_{i}[t]+
\frac{\beta[t+1]}{\sum_{r=1}^{t+1}\beta[r]}(\mathbf{\lambda}_{i}[t+1]
-\widehat{\mathbf{\lambda}}_{i}[t]),~\mathrm{for~all}~t\geq 0$$
where $\widehat{\mathbf{\lambda}}_{i}[0]=\lambda_{i}[0]$. Let us define the averaging iterates $\overline{\mu}[t]=\frac{\sum_{i=1}^{m}\mathbf{\mu}_{i}[t]}{m}$.\\
{\bf Remark 1.} i) Motivated by the algorithm proposed in \cite{nedic2015distributed}, we solve an optimization problem with coupling equality constraints by resorting to dual methods, while the problem considered in \cite{nedic2015distributed} has no coupling equality constraints. Our algorithm can be viewed as an extension of push-sum based algorithms \cite{nedic2015distributed} to a constrained setting.

ii) The primal-dual distributed methods proposed in \cite{zhu2012on,yuan2012distribted} require that each agent generates local copies of primal and dual variables, which then are optimized and exchanged. This, however, immediately leads to an increased computational and communication effort, which indeed scale as the number of agents.  In our method instead agents need to only optimize local variables and just exchange the estimate of dual variables, which are as many as the number of coupling constraints. The required local computational effort is thus much smaller when the number of coupling constraints is low compared to the overall dimensionality of primal decision variables.

iii) Under the assumption that the network is undirected, the work in \cite{falsone2016dual} obtain the convergence of dual optimality and primal optimality, without considering constraint violations of the problem interest. In our algorithm, we extend the algorithm in \cite{falsone2016dual} to directed graphs. We establish the convergence results for dual optimality, primal optimality and constraint violation, stated as in Theorem \ref{th1} below. More importantly, we derive the explicit convergence rate of the proposed algorithm, presented in the following Theorems \ref{th2} and \ref{th3}.

\subsection{Main results}

In this section, we show that the convergence results of the proposed Algorithm \ref{alg1}. Moreover, we provide the explicit estimates on the convergence rate of objective function' values and constraint violations.

The following Theorem \ref{th1} shows the optimality of  dual and primal variables under certain stepsize rule. Specifically, each local estimates $\mathbf{\lambda}_{i}[t]$ approaches the optimal dual vector $\mathbf{\lambda}^{*}$, while the auxiliary sequence $\{\widehat{\mathbf{x}}[t]=(\widehat{\mathbf{x}}_{1}[t]^{\top},\ldots,\widehat{\mathbf{x}}_{m}[t]^{\top})^{\top}\}$ converges to the optimal primal solution $\mathbf{x}^{*}$. In addition, the sequence $\{\widehat{\mathbf{x}}[t]\}$ satisfies the coupling equality constraints as $t\rightarrow\infty$.

\begin{theorem}\label{th1}
Suppose that Assumptions \ref{A1}-\ref{A3} hold, and the positive stepsize sequence $\{\beta[t]\}_{t\geq 1}$ satisfies the decay conditions below
\begin{eqnarray}\label{step}
\sum_{t=1}^{\infty}\beta[t]=\infty,~~\sum_{t=1}^{\infty}\beta^{2}[t]<\infty, ~~\beta[t]\leq\beta[l], ~~\textrm{for all} ~~t> l\geq 1.
\end{eqnarray}
Then, for some $\mathbf{x}^{*}\in \mathbf{X}^{*}$ and $\mathbf{\lambda}^{*}\in\Lambda^{*}$, we have that\\
(i) (dual optimality) $\lim_{t\rightarrow\infty}||\mathbf{\lambda}_{i}[t]-\mathbf{\lambda}^{*}||=0$, for all $i=1,2,\ldots,m;$\\
(ii) (feasibility) $\widehat{\mathbf{x}}[t]\in \mathbf{X}$ and $\lim_{t\rightarrow\infty}\sum_{i=1}^{m}A_{i}\widehat{\mathbf{x}}_{i}[t]-\mathbf{b}_{i}=0;$\\
(iii) (primal optimality) $\lim_{t\rightarrow\infty}||\widehat{\mathbf{x}}[t]-\mathbf{x}^{*}||=0,~\forall~ \mathbf{x}^{*}\in \mathbf{X}^{*}.$
\end{theorem}

Next Theorems \ref{th2} and \ref{th3} give the convergence rate of Algorithm \ref{alg1} under suitable choice of  stepsize, which characterizes the convergent speedup for the values of primal objective function and violations of coupling constraints, respectively.

\begin{theorem}\label{th2} (Convergence rate) Consider Assumptions \ref{A1}-\ref{A3} and let $\beta[t]={c}/{\sqrt{t}}$. Then, for any $t\geq1$ and $\mathbf{x}^{*}\in \mathbf{X}^{*}$, we have
\begin{eqnarray}
&&F(\widehat{\mathbf{x}}[t+1])-F(\mathbf{x}^{*})\nonumber\\
\leq&&\frac{m||\overline{\mu}[0]||_{1}}{2c\sqrt{t+1}}+\frac{c G^{2}(1+\ln(t+1))}{2m\sqrt{t+1}}\nonumber\\
&&+ \frac{16G\sum_{j=1}^{m}||\mu_{j}[0]||_{1}}{\xi(1-\eta)\sqrt{t+1}}+ \frac{16cpG^{2}(1+\ln t)}{\xi(1-\eta)\sqrt{t+1}},\nonumber
\end{eqnarray}
where $c>0$ is a constant, and $\xi>0$ and $\eta\in (0,1)$ satisfy $\xi\geq\frac{1}{m^{mB}}$, $\eta\leq(1-\frac{1}{m^{mB}})^\frac{1}{mB}$.
\end{theorem}
{\bf Remark 2.} i) Theorem \ref{th2} implies that the iterative sequence of network objective function $\{F(\widehat{\mathbf{x}}[t+1])\}$ converges to the optimal objective value $F(\mathbf{x}^{*})$, i.e.,
$\lim_{t\rightarrow\infty}F(\widehat{\mathbf{x}}[t+1]) = F(\mathbf{x}^{*}).$

ii) More importantly, Theorem \ref{th2} shows that the diffidence between the iterative sequence of primal objective function $\{F(\widehat{\mathbf{x}}[t+1])\}$ and the optimal value $F(\mathbf{x}^{*})$ converges at a rate of $O({\ln t}/{\sqrt{t}})$, that is
$$F(\widehat{\mathbf{x}}[t+1]) - F(\mathbf{x}^{*})=O\left(\frac{\ln t}{\sqrt{t}}\right)$$
with the constant depending on the initial values $\overline{\mu}[0]$ at the agents, the subgradient norm $G$ of dual function, and on both the speed $\eta$ of the network information diffusion and the imbalances $\xi$ of influence among the agents.

\begin{theorem}\label{th3} (Constraint violations) Consider Assumptions \ref{A1}-\ref{A3} and let $\beta[t]={c}/{\sqrt{t}}$. Then, for any $t\geq1$, we have
\begin{eqnarray}
&& ||\sum_{i=1}^{m}A_{i}\widehat{\mathbf{x}}_{i}[t+1]-\mathbf{b}_{i}||^{2}    \nonumber\\
\leq && \frac{4m^{2}||\overline{\mathbf{\mu}}[0]||_{1}}{c^{2}(t+1)}+\frac{2G^{2}(1+\ln(t+1))}{{t+1}}\nonumber\\
&&+ \frac{64Gm\sum_{j=1}^{m}||\mu_{j}[0]||_{1}}{c\xi(1-\eta)(t+1)}+\frac{64mpG^{2}(1+\ln t)}{\xi(1-\eta)(t+1)}.\nonumber
\end{eqnarray}
\end{theorem}

Theorem \ref{th3} provides that the constraint violation measured by $||\sum_{i=1}^{m}A_{i}\widehat{\mathbf{x}}_{i}[t+1]-\mathbf{b}_{i}||$ is also of
the order $O({\sqrt{\ln t}}/{\sqrt{t}})$.

\section{Proof of main results }\label{sec4}

We first prove the result that the dual optimal solutions are restricted in some specific sets, which will be useful to deduce the convergence of Algorithm \ref{alg1}.

\begin{lemma}\label{lemma1}
Under the Assumptions \ref{A1} and  \ref{A2}, the dual optimal solution $\mathbf{\lambda}^{*}$ of dual problem (\ref{dual}) is bounded, i.e., there an exist constant $C>0$ such that $0\leq ||\lambda^{*}||\leq C$.
\end{lemma}
{\bf Proof.} Letting $\check{\mathbf{x}}\in \mathbf{X}$ be a Slater vector, it holds that $\sum_{i}^{m}A_{i}\check{\mathbf{x}}_{i}-\mathbf{b}_{i}=\mathbf{0}$. Under Assumption \ref{A2}, the strong duality holds, that is, for all $\mathbf{\lambda} \in \mathbb{R}^{p}, \mathbf{x}\in \mathbf{X}$
\begin{eqnarray}\label{eq-le1-1}
L(\mathbf{x},\mathbf{\lambda}^{*})\geq  L(\mathbf{x}^{*},\mathbf{\lambda}^{*})\geq L(\mathbf{x}^{*},\mathbf{\lambda})\geq L(\mathbf{x}(\mathbf{\lambda}),\mathbf{\lambda})=\phi(\mathbf{\lambda}),
\end{eqnarray}
where $\mathbf{x}(\mathbf{\lambda})=\mathrm{arg}\min_{x\in \mathbf{X}} L(\mathbf{x},\mathbf{\lambda}) $.
By (\ref{eq-le1-1}), it gives rise to
\begin{eqnarray}\label{eq-le1-2}
\sum_{i=1}^{m}f_{i}(\mathbf{x}_{i}) + \mathbf{\lambda}^{*\top}(\sum_{i=1}^{m}A_{i}\mathbf{x}_{i}-\mathbf{b}_{i})\geq  \phi(\mathbf{\lambda}).\nonumber
\end{eqnarray}
Using the  inequality above and the fact $\sum_{i}^{m}\mathbf{b}_{i}= \sum_{i}^{m}A_{i}\check{\mathbf{x}}_{i}$, for any $ \mathbf{x}_{i}\in \mathbf{X}_{i}$ and $\mathbf{\lambda} \in \mathbb{R}^{p}$, we can obtain
\begin{eqnarray}\label{eq-le1-3}
\sum_{i=1}^{m}f_{i}(\mathbf{x}_{i}) + \mathbf{\lambda}^{*\top}\sum_{i=1}^{m}A_{i}(\mathbf{x}_{i}- \check{\mathbf{x}}_{i})\geq \phi(\mathbf{\lambda})=\sum_{i=1}^{m}\phi_{i}(\mathbf{\lambda}).
\end{eqnarray}

Letting $\mathbf{x}_{i}=\check{\mathbf{x}}_{i} + \mathbf{e}_{l_{i}}$ in (\ref{eq-le1-3}), where $\mathbf{e}_{l_{i}}$ denote the unit vector such that the $l_{i}$th component of $\mathbf{e}_{l_{i}}$ equals to 1, and the other components of $\mathbf{e}_{l_{i}}$ are 0, $l_{i}=1,2,\ldots, n_{i}$, we have
\begin{eqnarray}\label{eq-le1-4}
\mathbf{\lambda}^{*\top}\sum_{i=1}^{m}(A_{i}\mathbf{e}_{l_{i}}) \geq  \sum_{i=1}^{m}\phi_{i}(\mathbf{\lambda})-f_{i}(\check{\mathbf{x}}_{i}+\mathbf{e}_{l_{i}}),\nonumber
\end{eqnarray}
so
\begin{eqnarray}\label{eq-le1-5}
\mathbf{\lambda}^{*\top}\sum_{i=1}^{m}(A_{i}\mathbf{e}_{l_{i}}) \geq  \sum_{i=1}^{m} \min_{1\leq l_{i}\leq n_{i}}\{\phi_{i}(\mathbf{\lambda})-f_{i}(\check{\mathbf{x}}_{i}+\mathbf{e}_{l_{i}})\}.
\end{eqnarray}
Similarly, selecting $\mathbf{x}_{i}=\check{\mathbf{x}}_{i} - \mathbf{e}_{l_{i}}$, we obtain
\begin{eqnarray}\label{eq-le1-6}
\mathbf{\lambda}^{*\top}\sum_{i=1}^{m}(A_{i}\mathbf{e}_{l_{i}}) \leq \sum_{i=1}^{m} f_{i}(\check{\mathbf{x}}_{i}-\mathbf{e}_{l_{i}})-\phi_{i}(\mathbf{\lambda}), \nonumber
\end{eqnarray}
so
\begin{eqnarray}\label{eq-le1-7}
\mathbf{\lambda}^{*\top}\sum_{i=1}^{m}(A_{i}\mathbf{e}_{l_{i}}) \leq \sum_{i=1}^{m} \max_{1\leq l_{i}\leq n_{i}} \{f_{i}(\check{\mathbf{x}}_{i}-\mathbf{e}_{l_{i}})-\phi_{i}(\mathbf{\lambda})\}.
\end{eqnarray}
By (\ref{eq-le1-5}) and (\ref{eq-le1-7}), for all $ \lambda\in R^{p}$, we get
\begin{eqnarray}\label{eq-le1-8}
|\mathbf{\lambda}^{*\top}\sum_{i=1}^{m}(A_{i}\mathbf{e}_{l_{i}})| \leq \sum_{i=1}^{m} \max_{1\leq l_{i}\leq n_{i}}\{\max\{ f_{i}(\check{\mathbf{x}}_{i}-\mathbf{e}_{l_{i}})-\phi_{i}(\mathbf{\lambda}), f_{i}(\check{\mathbf{x}}_{i}+\mathbf{e}_{l_{i}})-\phi_{i}(\mathbf{\lambda})\}\}.\nonumber
\end{eqnarray}
Choosing a bounded vector $ \check{\mathbf{\lambda}}\in \mathbb{R}^{p}$ randomly and letting $\lambda=\check{\mathbf{\lambda}}$ in the above inequality, it leads to
\begin{eqnarray}\label{eq-le1-9}
|\mathbf{\lambda}^{*\top}\sum_{i=1}^{m}(A_{i}\mathbf{e}_{l_{i}})| \leq \sum_{i=1}^{m} \max_{1\leq l_{i}\leq n_{i}}\{\max\{ f_{i}(\check{\mathbf{x}}_{i}-\mathbf{e}_{l_{i}})-\phi_{i}(\check{\mathbf{\lambda}}), f_{i}(\check{\mathbf{x}}_{i}+\mathbf{e}_{l_{i}})-\phi_{i}(\check{\mathbf{\lambda}})\}\},
\end{eqnarray}
Letting $\breve{C}_{l_{i}}$ be an arbitrary value larger than $\max\{ f_{i}(\check{\mathbf{x}}_{i}-\mathbf{e}_{l_{i}})-\phi_{i}(\check{\mathbf{\lambda}}), f_{i}(\check{\mathbf{x}}_{i}+\mathbf{e}_{l_{i}})-\phi_{i}(\check{\mathbf{\lambda}})\}$, it follows from (\ref{eq-le1-9}) that
\begin{eqnarray}\label{eq-le1-10}
|\mathbf{\lambda}^{*\top}\sum_{i=1}^{m}(A_{i}\mathbf{e}_{l_{i}})|\leq \sum_{i=1}^{m}\max_{l_{i}}\breve{C}_{l_{i}}=\breve{C}.
\end{eqnarray}
Since the vector $\sum_{i=1}^{m}(A_{i}\mathbf{e}_{l_{i}})$ is a constant in (\ref{eq-le1-10}), $||\mathbf{\lambda}^{*}||$ is bounded.

Note that $\mathbf{x}_{i}=\check{\mathbf{x}}_{i} \pm \mathbf{e}_{l_{i}}$ may not belong to $\mathbf{X}_{i}$, but $\check{\mathbf{x}}_{i}$ is an
interior point of $\mathbf{X}_{i}$, so there exists a small number $\epsilon>0$ such that $\mathbf{x}_{i}=\check{\mathbf{x}}_{i} \pm \epsilon\mathbf{e}_{l_{i}}\in \mathbf{X}_{i}$. Then we can still have the same conclusion as above. The proof of this lemma  is completed.   \qed

Next we establish a fundamental lemma, which is helpful to prove the main results.

\begin{lemma}\label{lemma2}
Under the Assumptions \ref{A1}-\ref{A3}, for all $\mathbf{x}\in \mathbf{X}$ and $\mathbf{\lambda} \in \mathbb{R}^{p}$, we have,
\begin{eqnarray}
&&||\overline{\mathbf{\mu}}[t+1]-\mathbf{\lambda}||^{2}\nonumber\\
\leq&& ||\overline{\mathbf{\mu}}[t]-\mathbf{\lambda}||^{2}+
\frac{4\beta[t+1]}{m}\sum_{i=1}^{m}G_{i}||\mathbf{\lambda}_{i}[t+1]- \overline{\mathbf{\mu}}[t]|| \nonumber\\
&&+\frac{G^{2}}{m^{2}}\beta^{2}[t+1] - \frac{2\beta[t+1]}{m}(L(\mathbf{x}[t+1],\mathbf{\lambda})- L(\mathbf{x},\overline{\mathbf{\mu}}[t])).\nonumber
\end{eqnarray}
\end{lemma}
{\bf Proof.}  By Step 8 of Algorithm \ref{alg1} and the column-stochasticity of matrix $D[t]$ defined by (\ref{cop2}),  we can obtain
\begin{eqnarray}\label{eq-theorem4-1}
\overline{\mathbf{\mu}}[t+1]=\overline{\mathbf{\mu}}[t]+\frac{\beta[t+1]}{m}\sum_{j=1}^{m}(A_{j}\mathbf{x}_{j}[t+1]-\mathbf{b}_{j}).
\end{eqnarray}
Due to $||A_{i}\mathbf{x}_{i}-\mathbf{b}_{i}||\leq G_{i}$, for any $\mathbf{\lambda}\in \mathbb{R}^{p}$, it follows from (\ref{eq-theorem4-1}) that
\begin{eqnarray}\label{eq-theorem4-3}
||\overline{\mathbf{\mu}}[t+1]-\mathbf{\lambda}||^{2}&\leq& ||\overline{\mathbf{\mu}}[t]-\mathbf{\lambda}||^{2}
+ \frac{\beta^{2}[t+1]}{m^{2}}G^{2}\nonumber\\
&+& \frac{2\beta[t+1]}{m}\sum_{j=1}^{m}(A_{j}\mathbf{x}_{j}[t+1]-\mathbf{b}_{j})^{\top}(\overline{\mathbf{\mu}}[t]-\mathbf{\lambda}).
\end{eqnarray}
Considering the cross-term $(A_{j}\mathbf{x}_{j}[t+1]-\mathbf{b}_{j})^{\top}(\overline{\mathbf{\mu}}[t]-\mathbf{\lambda})$ in (\ref{eq-theorem4-3}), we have
\begin{eqnarray}\label{eq-theorem4-4}
(A_{j}\mathbf{x}_{j}[t+1]-\mathbf{b}_{j})^{\top}(\overline{\mathbf{\mu}}[t]-\mathbf{\lambda})&=&(A_{j}\mathbf{x}_{j}[t+1] -\mathbf{b}_{j})^{\top}(\overline{\mathbf{\mu}}[t]-
\mathbf{\lambda}_{j}[t+1])\nonumber\\
&+&(A_{j}\mathbf{x}_{j}[t+1]-\mathbf{b}_{j})^{\top}(\mathbf{\lambda}_{j}[t+1]-\mathbf{\lambda}).
\end{eqnarray}
Using the Cauchy-Schwartz inequality, we can get
\begin{eqnarray}\label{eq-theorem4-5}
(A_{j}\mathbf{x}_{j}[t+1]-\mathbf{b}_{j})^{\top}(\overline{\mathbf{\mu}}[t]-\mathbf{\lambda}_{j}[t+1])\leq G_{j}||\overline{\mathbf{\mu}}[t]-\mathbf{\lambda}_{j}[t+1]||.
\end{eqnarray}
For the second term of right-hand side in (\ref{eq-theorem4-4}), we can obtain
\begin{eqnarray}\label{eq-theorem4-6}
(A_{j}\mathbf{x}_{j}[t+1]-\mathbf{b}_{j})^{\top}(\mathbf{\lambda}_{j}[t+1]-\mathbf{\lambda})&=&f_{j}(\mathbf{x}_{j}[t+1])+
\mathbf{\lambda}_{j}[t+1]^{\top}(A_{j}\mathbf{x}_{j}[t+1]-\mathbf{b}_{j})\nonumber\\
&-&(f_{j}(\mathbf{x}_{j}[t+1])+\mathbf{\lambda}^{\top}(A_{j}\mathbf{x}_{j}[t+1]-\mathbf{b}_{j})).
\end{eqnarray}
By Step 7 of Algorithm \ref{alg1}, we get
$$f_{j}(\mathbf{x}_{j}[t+1])+\mathbf{\lambda}_{j}[t+1]^{\top}(A_{j}\mathbf{x}_{j}[t+1]-\mathbf{b}_{j})\leq f_{j}(\mathbf{x}_{j})+\mathbf{\lambda}_{j}[t+1]^{\top}(A_{j}\mathbf{x}_{j}-\mathbf{b}_{j}).$$
Using the inequality as above and (\ref{eq-theorem4-6}), we have
\begin{eqnarray}\label{eq-theorem4-7}
&&(A_{j}\mathbf{x}_{j}[t+1]-\mathbf{b}_{j})^{\top}(\mathbf{\lambda}_{j}[t+1]-\mathbf{\lambda})\nonumber\\
\leq&& f_{j}(\mathbf{x}_{j})+\mathbf{\lambda}_{j}[t+1]^{\top}(A_{j}\mathbf{x}_{j}-\mathbf{b}_{j})
-(f_{j}(\mathbf{x}_{j}[t+1])+\mathbf{\lambda}^{\top}(A_{j}\mathbf{x}_{j}[t+1]-\mathbf{b}_{j}))\nonumber\\
=&&f_{j}(\mathbf{x}_{j})+\mathbf{\lambda}_{j}[t+1]^{\top}(A_{j}\mathbf{x}_{j}-\mathbf{b}_{j})
-(f_{j}(\mathbf{x}_{j})+\overline{\mathbf{\mu}}[t]^{\top}(A_{j}\mathbf{x}_{j}-\mathbf{b}_{j}))\nonumber\\
&&+(f_{j}(\mathbf{x}_{j})+\overline{\mathbf{\mu}}[t]^{\top}(A_{j}\mathbf{x}_{j}-\mathbf{b}_{j}))
-(f_{j}(\mathbf{x}_{j}[t+1])+\mathbf{\lambda}^{\top}(A_{j}\mathbf{x}_{j}[t+1]-\mathbf{b}_{j}))\nonumber\\
\leq&& G_{j}||\overline{\mathbf{\mu}}[t]-\mathbf{\lambda}_{j}[t+1]||+L_{j}(\mathbf{x}_{j},\overline{\mathbf{\mu}}[t])
-L_{j}(\mathbf{x}_{j}[t+1],\mathbf{\lambda}),
\end{eqnarray}
where the last inequality uses the definition of the function $L_{j}$ given by (\ref{lf}). Finally, combining (\ref{eq-theorem4-3}), (\ref{eq-theorem4-4}), (\ref{eq-theorem4-5}) and (\ref{eq-theorem4-7}), we can obtain the conclusion.   \qed

In what follows, we give the well-known Supermartingale Convergence Theorem \cite{polyak1987introduction}, refer to Lemma \ref{lemma3}, which is useful to prove Theorem \ref{th1}.

\begin{lemma}\label{lemma3} Let $\{x[t]\}$ be a non-negative scalar sequence such that
$$x[t+1]\leq(1+b[t])x[t]-y[t]+c[t],~\forall t\geq 0$$
where $b[t]\geq0$, $y[t]\geq0$ and $c[t]\geq0$ for all $t\geq0$ with $\sum_{t=0}^{\infty}b[t]<\infty$, and $\sum_{t=0}^{\infty}c[t]<\infty$. Then, the sequence $\{x[t]\}$ converges to some $x\geq0$ and $\sum_{t=0}^{\infty}y[t]<\infty$.
\end{lemma}

We are ready to give the proof of main results. We firstly prove Theorem \ref{th1} to show the optimality of dual and primal variables, and constraint violations.\\\\
{\bf Proof of Theorem \ref{th1}}: Letting $\mathbf{\lambda}=\mathbf{\lambda}^{*}$ and $\mathbf{x}=\mathbf{x}^{*}$ in Lemma \ref{lemma2}, for some $\mathbf{\lambda}^{*}\in \Lambda^{*}$ and $\mathbf{x}^{*}\in \mathbf{X}^{*}$, we have
\begin{eqnarray}\label{eq-theorem1-3}
&&||\overline{\mathbf{\mu}}[t+1]-\mathbf{\lambda}^{*}||^{2}\nonumber\\
\leq&&||\overline{\mathbf{\mu}}[t]-\mathbf{\lambda}^{*}||^{2}
+\frac{4\beta[t+1]}{m}\sum_{i=1}^{m}G_{i}||\mathbf{\lambda}_{i}[t+1]-\overline{\mathbf{\mu}}[t]||\nonumber\\
&&+ \frac{G^{2}}{m^{2}}\beta^{2}[t+1]
 -\frac{2\beta[t+1]}{m}(L(\mathbf{x}[t+1],\mathbf{\lambda}^{*})-L(\mathbf{x}^{*},\overline{\mathbf{\mu}}[t])).
\end{eqnarray}
Making use of the saddle-point theorem
and (\ref{eq-theorem1-3}), it gives rise to
\begin{eqnarray}\label{eq-theorem1-4}
&&||\overline{\mathbf{\mu}}[t+1]-\mathbf{\lambda}^{*}||^{2}\nonumber\\
\leq&&||\overline{\mathbf{\mu}}[t]-\mathbf{\lambda}^{*}||^{2}
+\frac{4\beta[t+1]}{m}\sum_{i=1}^{m}G_{i}||\mathbf{\lambda}_{i}[t+1]-\overline{\mathbf{\mu}}[t]||  \nonumber\\
&&+ \frac{G^{2}}{m^{2}}\beta^{2}[t+1]
- \frac{2\beta[t+1]}{m}(L(\mathbf{x}^{*},\mathbf{\lambda}^{*})-L(\mathbf{x}^{*},\overline{\mathbf{\mu}}[t])),
\end{eqnarray}
According to Lemma 1 (b) in \cite{nedic2015distributed}, the following result holds
\begin{eqnarray}\label{eq-theorem1-1}
\lim_{t\rightarrow\infty}||\mathbf{\lambda}_{i}[t+1]-\overline{\mathbf{\mu}}[t]||=0, ~\forall~ i=1,2,\ldots,m.
\end{eqnarray}
Since $\sum_{t=1}^{\infty}\beta[t]<\infty$, it follows from (\ref{eq-theorem1-1}) that
\begin{eqnarray}
\sum_{t=0}^{\infty}\beta[t+1]||\mathbf{\lambda}_{i}[t+1]-\overline{\mathbf{\mu}}[t]||<\infty, ~\forall~ i=1,2,\ldots,m. \nonumber
\end{eqnarray}
Further, we have
\begin{eqnarray}\label{eq-theorem1-2}
\sum_{t=0}^{\infty}\frac{4\beta[t+1]}{m}\sum_{i=1}^{m}G_{i}||\mathbf{\lambda}_{i}[t+1]-\overline{\mathbf{\mu}}[t]||<\infty.
\end{eqnarray}
Note that the fact that $L(\mathbf{x}^{*},\mathbf{\lambda}^{*})-L(\mathbf{x}^{*},\overline{\mathbf{\mu}}[t])\geq 0$. Thus, by Lemma \ref{lemma3},  we can conclude that the sequence $\{\overline{\mathbf{\mu}}[t]\}$ converges to the solution $\mathbf{\lambda}^{*}\in \Lambda^{*}$. Furthermore, by (\ref{eq-theorem1-1}), we can see that each sequence $\{\lambda_{i}[t]\}$ converges to the same solution $\mathbf{\lambda}^{*}$, $i=1,2,\ldots,m$, thus concluding the proof of i) in Theorem \ref{th1}.

Next we prove the ii) of Theorem \ref{th1}.
By the definition of $\widehat{\mathbf{x}}_{i}[t+1]$, it can be rewritten as
\begin{eqnarray}\label{eq-theorem1-5}
\widehat{\mathbf{x}}_{i}[t+1]
=\frac{\sum_{r=0}^{t}\beta[r+1]\mathbf{x}_{i}[r+1]}{\sum_{r=0}^{t}\beta[r+1]},
\end{eqnarray}
implying that $\widehat{\mathbf{x}}_{i}[t+1]$ is a convex combination of past values of $\mathbf{x}_{i}[t+1]$. Thus, for all $t\geq0$, we have that $\widehat{\mathbf{x}}_{i}[t+1]\in \mathbf{X}_{i}$ and
\begin{eqnarray}\label{eq-theorem1-6}
\sum_{i=1}^{m}(A_{i}\widehat{\mathbf{x}}_{i}[t+1]-\mathbf{b}_{i})&=&\sum_{i=1}^{m}[A_{i}(\frac{\sum_{r=0}^{t}\beta[r+1]\mathbf{x}_{i}[r+1]}
{\sum_{r=0}^{t}\beta[r+1]})-\mathbf{b}_{i}]\nonumber\\
&=&\frac{\sum_{i=1}^{m}\sum_{r=0}^{t}\beta[r+1](A_{i}\mathbf{x}_{i}[r+1]-\mathbf{b}_{i})}{\sum_{r=0}^{t}\beta[r+1]}.\nonumber
\end{eqnarray}
By Step 8 of Algorithm \ref{alg1} and the column-stochasticity of matrix $D[t]$, it follows that
\begin{eqnarray}\label{eq-theorem1-7}
\sum_{i=1}^{m}(A_{i}\widehat{\mathbf{x}}_{i}[t+1]-\mathbf{b}_{i})&=&
\frac{\sum_{i=1}^{m}\sum_{r=0}^{t}(\mathbf{\mu}_{i}[r+1]-\mathbf{u}_{i}[r+1])}{\sum_{r=0}^{t}\beta[r+1]}\nonumber\\
&=& \frac{\sum_{i=1}^{m}\sum_{r=0}^{t}(\mathbf{\mu}_{i}[r+1]- \sum_{j=1}^{m}(D[r])_{ij}\mathbf{\mu}_{j}[r])}{\sum_{r=0}^{t}\beta[r+1]}\nonumber\\
&=& \frac{\sum_{r=0}^{t}(\sum_{i=1}^{m}\mathbf{\mu}_{i}[r+1]- \sum_{i=1}^{m}\mathbf{\mu}_{i}[r])}{\sum_{r=0}^{t}\beta[r+1]}\nonumber\\
&=& \frac{m(\overline{\mathbf{\mu}}[t+1]- \overline{\mathbf{\mu}}[0])}{\sum_{r=0}^{t}\beta[r+1]}.
\end{eqnarray}
Since $\lim_{t\rightarrow\infty}\overline{\mathbf{\mu}}[t+1]=\mathbf{\lambda}^{*}$, we further get
$$\lim_{t\rightarrow \infty}m(\overline{\mathbf{\mu}}[t+1]- \overline{\mathbf{\mu}}[0])< \infty.$$
Using $\sum_{r=0}^{\infty}\beta[r]=\infty$, it holds that
\begin{eqnarray}\label{eq-theorem1-8}
\lim_{t\rightarrow\infty}\frac{m(\overline{\mathbf{\mu}}[t+1]- \overline{\mathbf{\mu}}[0])}{\sum_{r=0}^{t}\beta[r+1]}=0.\nonumber
\end{eqnarray}
It follows from (\ref{eq-theorem1-7}) and the above relation that
\begin{eqnarray}\label{eq-theorem1-9}
\lim_{t\rightarrow\infty}\sum_{i=1}^{m}(A_{i}\widehat{\mathbf{x}}_{i}[t+1]-\mathbf{b}_{i})=0. \nonumber
\end{eqnarray}
which completes the proof of ii) in Theorem \ref{th1}.

Now we begin to prove the iii) of Theorem \ref{th1}.
Considering the quantity $2\sum_{i=1}^{m}L_{i}(\widehat{\mathbf{x}}_{i}[t+1],\mathbf{\lambda}^{*})$, and using the convexity of $L(\cdot,\mathbf{\lambda})$ and (\ref{eq-theorem1-5}), we obtain
\begin{eqnarray}\label{eq-theorem1-10}
2\sum_{i=1}^{m}L_{i}(\widehat{\mathbf{x}}_{i}[t+1],\mathbf{\lambda}^{*})&\leq&
\frac{2\sum_{i=1}^{m}\sum_{r=0}^{t}\beta[r+1]L_{i}(\mathbf{x}_{i}[r+1],\mathbf{\lambda}^{*})}{\sum_{r=0}^{t}\beta[r+1]}\nonumber\\
&=& \frac{\sum_{r=0}^{t}2\beta[r+1]L(\mathbf{x}[r+1],\mathbf{\lambda}^{*})}{\sum_{r=0}^{t}\beta[r+1]}.
\end{eqnarray}
Rearranging the terms in (\ref{eq-theorem1-3}) and  letting $t=r$, we have
\begin{eqnarray}\label{eq-theorem1-12}
&&2\beta[r+1]L(\mathbf{x}[r+1],\mathbf{\lambda}^{*})\nonumber\\
\leq &&2\beta[r+1]L(\mathbf{x}^{*},\overline{\mathbf{\mu}}[r])+\frac{G^{2}}{m}\beta^{2}[r+1]+ m(||\overline{\mathbf{\mu}}[r]-\mathbf{\lambda}^{*}||^{2}-||\overline{\mathbf{\mu}}[r+1]-\mathbf{\lambda}^{*}||^{2}) \nonumber\\
&&
+4\beta[r+1]\sum_{i=1}^{m}G_{i}||\mathbf{\lambda}_{i}[r+1]-\overline{\mathbf{\mu}}[r]||\nonumber\\
\leq&&2\beta[r+1]L(\mathbf{x}^{*},\mathbf{\lambda}^{*})+\frac{G^{2}}{m}\beta^{2}[r+1]
+ m(||\overline{\mathbf{\mu}}[r]-\mathbf{\lambda}^{*}||^{2}-||\overline{\mathbf{\mu}}[r+1]-\mathbf{\lambda}^{*}||^{2})\nonumber\\
&&+ 4\beta[r+1]\sum_{i=1}^{m}G_{i}||\mathbf{\lambda}_{i}[r+1]-\overline{\mathbf{\mu}}[r]||,
\end{eqnarray}
where the second inequality is due to $L(\mathbf{x}^{*},\overline{\mathbf{\mu}}[r])\leq L(\mathbf{x}^{*},\mathbf{\lambda}^{*})$. Combining  (\ref{eq-theorem1-12})
and (\ref{eq-theorem1-10}), we can obtain
\begin{eqnarray}\label{eq-theorem1-13}
&&2L(\widehat{\mathbf{x}}[t+1],\mathbf{\lambda}^{*}) 2L(\mathbf{x}^{*},\mathbf{\lambda}^{*})\leq
\frac{\sum_{r=0}^{t}4\beta[r+1]\sum_{i=1}^{m}G_{i}||\mathbf{\lambda}_{i}[r+1]-\overline{\mathbf{\mu}}[r]||}{\sum_{r=0}^{t}\beta[r+1]}\nonumber\\
&&+ \frac{1}{\sum_{r=0}^{t}\beta[r+1]}\sum_{r=0}^{t}\frac{G^{2}}{m}\beta^{2}[r+1]
+\frac{m(||\overline{\mathbf{\mu}}[0]-\mathbf{\lambda}^{*}||^{2}-||\overline{\mathbf{\mu}}[t+1]-\mathbf{\lambda}^{*}||^{2}) }{\sum_{r=0}^{t}\beta[r+1]}.\nonumber\\
\end{eqnarray}
By (\ref{eq-theorem1-2}), the stepsize rule (\ref{step}) and the fact that  $\lim_{t\rightarrow\infty}\overline{\mathbf{\mu}}[t+1]=\mathbf{\lambda}^{*}$,
each term of right-hand side in (\ref{eq-theorem1-13}) is convergent to zero as $t\rightarrow\infty$, thus, we can obtain
\begin{eqnarray}\label{eq-theorem1-15}
\lim_{t\rightarrow\infty}\sup L(\widehat{\mathbf{x}}[t+1],\mathbf{\lambda}^{*})\leq L(\mathbf{x}^{*},\mathbf{\lambda}^{*}).\nonumber
\end{eqnarray}
Note that $L(\widehat{\mathbf{x}}[t+1],\mathbf{\lambda}^{*})\geq L(\mathbf{x}^{*},\mathbf{\lambda}^{*}) $, we further get
\begin{eqnarray}
\lim_{t\rightarrow\infty}\sup L(\widehat{\mathbf{x}}[t+1],\mathbf{\lambda}^{*})= L(\mathbf{x}^{*},\mathbf{\lambda}^{*})=F(\mathbf{x}^{*}).\nonumber
\end{eqnarray}
Since $ L(\cdot,\mathbf{\lambda}^{*})$ is continuous and convex for any $\mathbf{\lambda}^{*}$,  all limit points of $\{\widehat{\mathbf{x}}[t]\}_{t\rightarrow\infty}$ are feasible and achieve the optimal value. This means that these limit points are optimal for the primal problem, thus, the proof of Theorem \ref{th1} is completed. \qed

Theorem \ref{th2} shows the convergence rate of the objective function's value under Assumptions \ref{A1}-\ref{A3}.\\\\
{\bf Proof of Theorem \ref{th2}}: By Lemma \ref{lemma2}, for any $\mathbf{x}\in \mathbf{X}$ and $\mathbf{\lambda}\in \mathbb{R}^{p}$, we can obtain
\begin{eqnarray}\label{eq-theorem2-1}
&&\sum_{r=0}^{t}\frac{2\beta[r+1]}{m}(L(\mathbf{x}[r+1],\mathbf{\lambda})-L(\mathbf{x},\overline{\mathbf{\mu}}[r]))\nonumber\\
\leq&& ||\overline{\mathbf{\mu}}[0]-\mathbf{\lambda}||^{2}-||\overline{\mathbf{\mu}}[t+1]-\mathbf{\lambda}||^{2}
+ \sum_{r=0}^{t}\frac{4\beta[r+1]}{m}\sum_{i=1}^{m}G_{i}||\overline{\mathbf{\mu}}[r]-\mathbf{\lambda}_{i}[r+1]||\nonumber\\
&&+ \frac{G^{2}}{m^{2}}\sum_{r=0}^{t}\beta^{2}[r+1]\nonumber\\
\leq&& ||\overline{\mathbf{\mu}}[0]-\mathbf{\lambda}||^{2}+ \sum_{r=0}^{t}\frac{4\beta[r+1]}{m}\sum_{i=1}^{m}G_{i}||\overline{\mathbf{\mu}}[r]-\mathbf{\lambda}_{i}[r+1]||
+\frac{G^{2}}{m^{2}}\sum_{r=0}^{t}\beta^{2}[r+1].\nonumber
\end{eqnarray}
Dividing both sides in above inequality by $\frac{2}{m}\sum_{r=0}^{t}\beta[r+1]$,  for any $\mathbf{x}\in \mathbf{X}$ , $\mathbf{\lambda}\in \mathbb{R}^{p}$ and $t>0$, we have
\begin{eqnarray}\label{eq-theorem2-2}
&&\frac{\sum_{r=0}^{t}\beta[r+1](L(\mathbf{x}[r+1],\mathbf{\lambda})-L(\mathbf{x},\overline{\mathbf{\mu}}[r]))}{\sum_{r=0}^{t}\beta[r+1]}\nonumber\\
\leq&& \frac{m||\overline{\mathbf{\mu}}[0]-\mathbf{\lambda}||^{2}}{2\sum_{r=0}^{t}\beta[r+1]}+\frac{1}{\sum_{r=0}^{t}\beta[r+1]}\sum_{r=0}^{t}\frac{G^{2}}{2m}\beta^{2}[r+1]\nonumber\\
&&+\frac{2}{\sum_{r=0}^{t}\beta[r+1]}\sum_{r=0}^{t}\beta[r+1]\sum_{i=1}^{m}G_{i}||\overline{\mu}[r]-\mathbf{\lambda}_{i}[r+1]||.
\end{eqnarray}
Note that  $L(\cdot,\mathbf{\lambda})$ is convex, it follows that
\begin{eqnarray}\label{eq-theorem2-4}
L(\frac{\sum_{r=0}^{t}\beta[r+1]\mathbf{x}[r+1]}{\sum_{r=0}^{t}\beta[r+1]},\mathbf{\lambda})
\leq\frac{\sum_{r=0}^{t}\beta[r+1]L(\mathbf{x}[r+1],\mathbf{\lambda})}{\sum_{r=0}^{t}\beta[r+1]}.\nonumber
\end{eqnarray}
Similarly, due to the concavity of $L(\mathbf{x},\cdot)$, it holds
\begin{eqnarray}\label{eq-theorem2-5}
\frac{\sum_{r=0}^{t}\beta[r+1]L(\mathbf{x},\overline{\mathbf{\mu}}[r])}{\sum_{r=0}^{t}\beta[r+1]}\leq L(\mathbf{x},\frac{\sum_{r=0}^{t}\beta[r+1]\overline{\mathbf{\mu}}[r]}{\sum_{r=0}^{t}\beta[r+1]}),\nonumber
\end{eqnarray}
Thus, for all $\mathbf{x}\in \mathbf{X}$, $\lambda \in R^{p}$, we obtain
\begin{eqnarray}
&&\frac{\sum_{r=0}^{t}\beta[r+1](L(\mathbf{x}[r+1],\mathbf{\lambda})-L(\mathbf{x},\overline{\mathbf{\mu}}[r]))}{\sum_{r=0}^{t}\beta[r+1]}\nonumber\\
&&\geq L(\frac{\sum_{r=0}^{t}\beta[r+1]\mathbf{x}[r+1]}{\sum_{r=0}^{t}\beta[r+1]},\mathbf{\lambda})-L(\mathbf{x},\frac{\sum_{r=0}^{t}\beta[r+1]\overline{\mathbf{\mu}}[r]}{\sum_{r=0}^{t}\beta[r+1]}).
\nonumber \end{eqnarray}
Letting $\mathbf{x}=\mathbf{x}^{*}$, $\mathbf{\lambda}=\mathbf{0}$ in the above inequality and using $L(\mathbf{x}^{*},\mathbf{\lambda})\leq L(\mathbf{x}^{*},\mathbf{\lambda}^{*})$, it gives rise to
\begin{eqnarray}\label{eq-theorem3-4}
&&\frac{\sum_{r=0}^{t}\beta[r+1](L(\mathbf{x}[r+1],\mathbf{0})-L(\mathbf{x}^{*},\overline{\mathbf{\mu}}[r]))}{\sum_{r=0}^{t}\beta[r+1]}\nonumber\\
&&\geq L(\frac{\sum_{r=0}^{t}\beta[r+1]\mathbf{x}[r+1]}{\sum_{r=0}^{t}\beta[r+1]},\mathbf{0})-
L(\mathbf{x}^{*},\frac{\sum_{r=0}^{t}\beta[r+1]\overline{\mathbf{\mu}}[r]}{\sum_{r=0}^{t}\beta[r+1]})\nonumber\\
&&\geq L(\frac{\sum_{r=0}^{t}\beta[r+1]\mathbf{x}[r+1]}{\sum_{r=0}^{t}\beta[r+1]},\mathbf{0})-L(\mathbf{x}^{*},\mathbf{\lambda}^{*})\nonumber\\
&&=F(\widehat{\mathbf{x}}[t+1])-F(\mathbf{x}^{*}).
\end{eqnarray}
Combining  (\ref{eq-theorem2-2}) and (\ref{eq-theorem3-4}), it yields
\begin{eqnarray}\label{eq-theorem2-6}
&&F(\widehat{\mathbf{x}}[t+1])-F(\mathbf{x}^{*})\nonumber\\
\leq&&\frac{m||\overline{\mathbf{\mu}}[0]||^{2}}{2\sum_{r=0}^{t}\beta[r+1]}+\frac{1}{\sum_{r=0}^{t}\beta[r+1]}\sum_{r=0}^{t}\frac{G^{2}}{2m}\beta^{2}[r+1]\nonumber\\
&&+\frac{2}{\sum_{r=0}^{t}\beta[r+1]}\sum_{r=0}^{t}\beta[r+1]\sum_{i=1}^{m}G_{i}||\overline{\mu}[r]-\mathbf{\lambda}_{i}[r+1]||.
\end{eqnarray}
Letting $\beta[t]=\frac{c}{\sqrt{t}}$, similar to Corollary 3 in \cite{nedic2015distributed}, we can deduce the following result
\begin{eqnarray}\label{eq-theorem2-7}
\sum_{r=0}^{t}\frac{c}{\sqrt{r+1}}||\overline{\mu}[r]-\mathbf{\lambda}_{i}[r+1]||
\leq \frac{8}{\xi(1-\eta)}(c||\mu_{j}[0]||_{1}+c^{2}p G(1+\ln t)).
\end{eqnarray}
Furthermore,
\begin{eqnarray}\label{eq-theorem2-8}
\sum_{r=0}^{t}\beta[r+1]=\sum_{r=0}^{t}\frac{c}{\sqrt{r+1}}\geq c\sqrt{t+1},
\end{eqnarray}
and
\begin{eqnarray}\label{eq-theorem2-9}
\sum_{r=0}^{t}\beta^{2}[t+1]=\sum_{l=1}^{t+1}\frac{c^{2}}{l}\leq c^{2}(1+\int_{1}^{t+1}\frac{dx}{x})=c^{2}(1+\ln (t+1)).
\end{eqnarray}
By (\ref{eq-theorem2-6}), (\ref{eq-theorem2-7}), (\ref{eq-theorem2-8}) and (\ref{eq-theorem2-9}), we have
\begin{eqnarray}\label{eq-theorem3-6}
&&F(\widehat{\mathbf{x}}[t+1])-F(\mathbf{x}^{*})\nonumber\\
\leq&& \frac{m||\overline{\mathbf{\mu}}[0]||_{1}}{2c\sqrt{t+1}}+\frac{G^{2}c(1+\ln(t+1))}{2m\sqrt{t+1}}\nonumber\\
&&+ \frac{16G\sum_{j=1}^{m}||\mu_{j}[0]||_{1}}{\xi(1-\eta)\sqrt{t+1}}+\frac{16pG^{2}c(1+\ln t)}{\xi(1-\eta)\sqrt{t+1}}.\nonumber
\end{eqnarray}
The proof of Theorem \ref{th2} is completed. \qed

Next, we begin to estimate the convergence rate of constraint violations.\\
{\bf Proof of Theorem \ref{th3}}: Letting $\lambda=0$ in (\ref{eq-theorem2-1}) and using (\ref{eq-theorem2-6}), we can deduce
\begin{eqnarray}
&&\frac{m||\overline{\mathbf{\mu}}[t+1]||^{2}}{2\sum_{r=0}^{t}\beta[r+1]}+F(\widehat{\mathbf{x}}[t+1])-F(\mathbf{x}^{*})\nonumber\\
\leq&&\frac{m||\overline{\mathbf{\mu}}[0]||^{2}}{2\sum_{r=0}^{t}\beta[r+1]}+\frac{1}{\sum_{r=0}^{t}\beta[r+1]}\sum_{r=0}^{t}\frac{G^{2}}{2m}\beta^{2}[r+1]\nonumber\\
&&+\frac{2}{\sum_{r=0}^{t}\beta[r+1]}\sum_{r=0}^{t}\beta[r+1]\sum_{i=1}^{m}G_{i}||\overline{\mu}[r]-\mathbf{\lambda}_{i}[r+1]||. \nonumber
\end{eqnarray}
Noting that $F(\widehat{\mathbf{x}}[t+1])-F(\mathbf{x}^{*})\geq0$, the above inequality leads to
\begin{eqnarray}\label{eq-theoremc-2}
\frac{m||\overline{\mathbf{\mu}}[t+1]||^{2}}{2\sum_{r=0}^{t}\beta[r+1]}
&\leq&\frac{m||\overline{\mathbf{\mu}}[0]||^{2}}{2\sum_{r=0}^{t}\beta[r+1]}+\frac{1}{\sum_{r=0}^{t}\beta[r+1]}\sum_{r=0}^{t}\frac{G^{2}}{2m}\beta^{2}[r+1]\nonumber\\
&&+\frac{2}{\sum_{r=0}^{t}\beta[r+1]}\sum_{r=0}^{t}\beta[r+1]\sum_{i=1}^{m}G_{i}||\overline{\mu}[r]-\mathbf{\lambda}_{i}[r+1]||.
\end{eqnarray}
Due to the fact that $||\overline{\mathbf{\mu}}[t+1]-\overline{\mathbf{\mu}}[0]||^{2} \leq 2(||\overline{\mathbf{\mu}}[t+1]||^{2}+||\overline{\mathbf{\mu}}[0]||^{2})$, it follows from (\ref{eq-theoremc-2}) that
\begin{eqnarray}\label{eq-theoremc-3}
&&\frac{m||\overline{\mathbf{\mu}}[t+1]-\overline{\mathbf{\mu}}[0]||^{2}}{\sum_{r=0}^{t}\beta[r+1]}\nonumber\\
\leq&&\frac{4m||\overline{\mathbf{\mu}}[0]||^{2}}{\sum_{r=0}^{t}\beta[r+1]}+\frac{2}{\sum_{r=0}^{t}\beta[r+1]}\sum_{r=0}^{t}\frac{G^{2}}{m}\beta^{2}[r+1]\nonumber\\
&&+\frac{8}{\sum_{r=0}^{t}\beta[r+1]}\sum_{r=0}^{t}\beta[r+1]\sum_{i=1}^{m}G_{i}||\overline{\mu}[r]-\mathbf{\lambda}_{i}[r+1]||.
\end{eqnarray}
Dividing both sides in (\ref{eq-theoremc-3}) by $\frac{\sum_{r=0}^{t}\beta[r+1]}{m}$, we have
\begin{eqnarray}\label{eq-theoremc-6}
&&\frac{m^{2}||\overline{\mathbf{\mu}}[t+1]-\overline{\mathbf{\mu}}[0]||^{2}}{(\sum_{r=0}^{t}\beta[r+1])^{2}}\nonumber\\
\leq&&\frac{4m^{2}||\overline{\mathbf{\mu}}[0]||^{2}}{(\sum_{r=0}^{t}\beta[r+1])^{2}}+\frac{2G^{2}}{(\sum_{r=0}^{t}\beta[r+1])^{2}}\sum_{r=0}^{t}\beta^{2}[r+1]\nonumber\\
&&+\frac{8m}{(\sum_{r=0}^{t}\beta[r+1])^{2}}\sum_{r=0}^{t}\beta[r+1]\sum_{i=1}^{m}G_{i}||\overline{\mu}[r]-\mathbf{\lambda}_{i}[r+1]||.
\end{eqnarray}
Combining (\ref{eq-theorem2-7}), (\ref{eq-theorem2-8}), (\ref{eq-theorem2-9}) and (\ref{eq-theoremc-6}), we can obtain
\begin{eqnarray}\label{eq-theoremc-5}
&&\frac{m^{2}||\overline{\mathbf{\mu}}[t+1]-\overline{\mathbf{\mu}}[0]||^{2}}{(\sum_{r=0}^{t}\beta[r+1])^{2}}\nonumber\\
\leq&& \frac{4m^{2}||\overline{\mathbf{\mu}}[0]||_{1}}{c^{2}(t+1)}+\frac{2G^{2}(1+\ln(t+1))}{{t+1}}\nonumber\\
&&+ \frac{64Gm\sum_{j=1}^{m}||\mu_{j}[0]||_{1}}{c\xi(1-\eta)(t+1)}+\frac{64mpG^{2}(1+\ln t)}{\xi(1-\eta)(t+1)}.
\end{eqnarray}
By (\ref{eq-theorem1-7}) and (\ref{eq-theoremc-5}), we can obtain the desired result.      \qed

\section{Numerical simulations}\label{sec5}

In this section we report and illustrate some experimental results of the proposed algorithm.

We consider the economic dispatch problem (EDP), which is vital in power system operation \cite{Xia2010}. EDP is commonly formulated as an optimization problem, where the objective is to minimize the total generation cost while meeting total demand and subjecting to individual generator output constraints:
\begin{eqnarray}\label{problem1}
  \min_{\{p_{i}\}_{i=1}^{m}} && C(p)=:\sum_{i=1}^{m}C_{i}(p_{i})\nonumber\\
  \textrm{s.t.} && \sum_{i=1}^{m} p_{i}=D,~~ p_{i}\in [p_{i}^{\mathrm{min}}, p_{i}^{\mathrm{max}}],~i=1,\ldots,m, \nonumber
\end{eqnarray}
where $p_{i}$ is the power generation of generator $i$, $m$ is the number of generators, $p_{i}^{\mathrm{min}}/p_{i}^{\mathrm{max}}$ are the corresponding minimum/maximum generation output, $D$
is the total demand, $C_{i}$ is the cost function of generator $i$.

In practice, individual generator $i$ only holds itself private information including the cost $C_{i}$ and generation limits $p_{i}^{\mathrm{min}},p_{i}^{\mathrm{max}}$ for privacy concern. Thus, centralized optimization methods are often unavailable.

The proposed Algorithm DDSG-PS is examined on IEEE 57-bus test system with 7 generators \cite{Han2001}. The cost function of each generator $i$ is taken as $C_{i}(p_{i})=a_{i}p^{2}_{i}+b_{i}p_{i}+c_{i}$, and the parameters of all the generators are given in Table \ref{table-1} \cite{Han2001}. The local demand at each generator is set as (241.0712,100.0000,74.8088,100.0000,550.0000,100.0000,410.0000) (MW) with total demand $D=1575.88$ MW.
\begin{table} [htb]
\begin{center} {\scriptsize}
\caption{Parameters for 7 generators in IEEE 57-bus system}\label{table-1}
\begin{tabular}{c|c|c|c|c}
  \hline \hline
Gen.  & $a_{i}$ & $b_{i}$  & $c_{i}$ & [$p_{i}^{\mathrm{min}},p_{i}^{\mathrm{max}}$]   \\   \hline
1  & 0.0775795		& 20  & 0  & [0,575.88] \\ \hline
2  & 0.01		    & 40  & 0  & [0,100] \\ \hline
3  & 0.25		    & 20  & 0  & [0,140] \\ \hline
6  & 0.01		    & 40  & 0  & [0,100] \\ \hline
8  & 0.0222222		& 20  & 0  & [0,550] \\ \hline
9  & 0.01		    & 40  & 0  & [0,100] \\ \hline
12  & 0.0322581		& 20  & 0  & [0,410] \\ \hline
\hline
\end{tabular}
\end{center}
\end{table}

\begin{figure}[!htbp]
\centering
\subfigure[Dual variables (\$/MWh) v.s. iterations ]{\includegraphics[width=18cm]{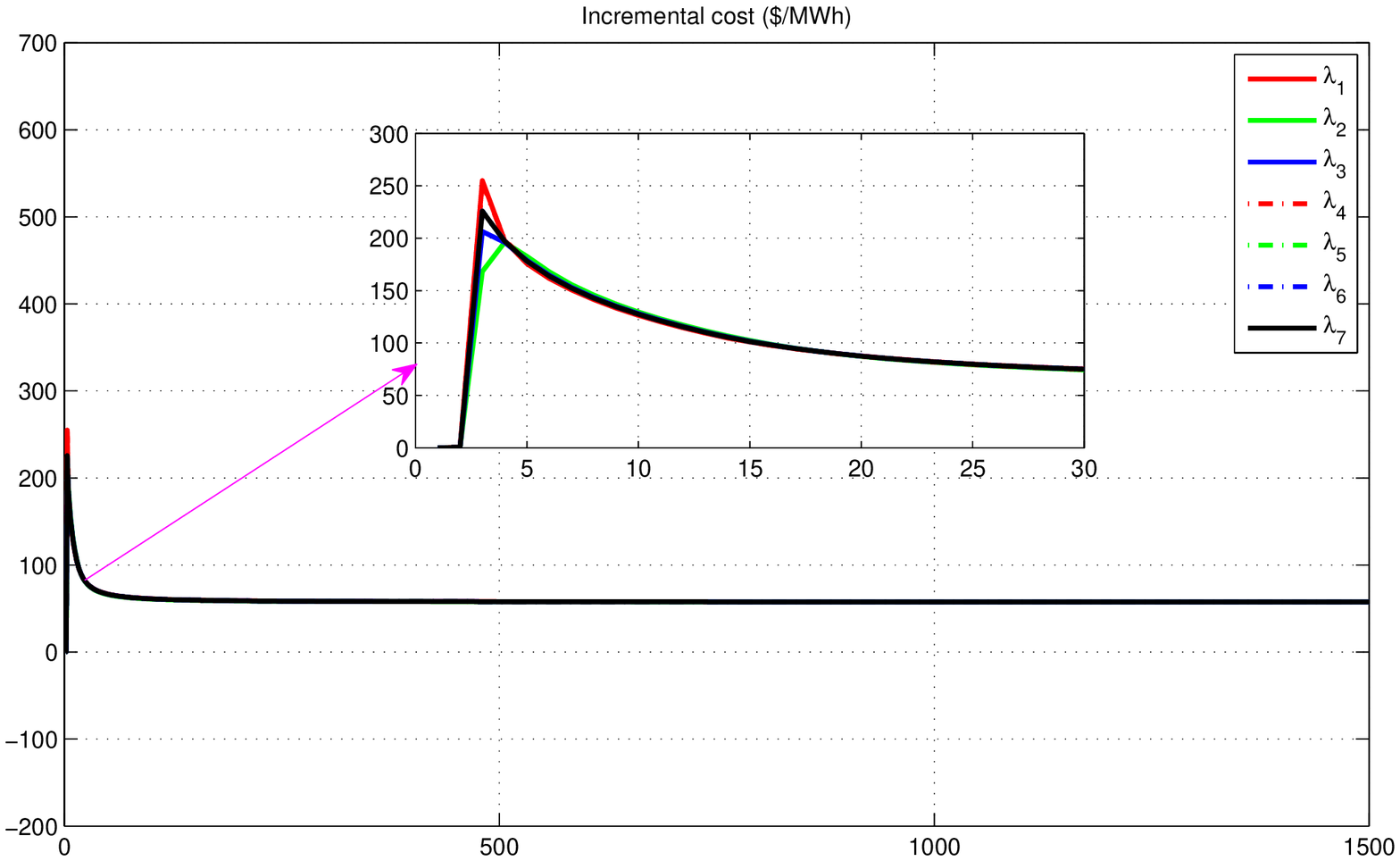}\label{fig1:subfig:a}}
\subfigure[Generation outputs (MW) v.s. iterations  ]{\includegraphics[width=18cm]{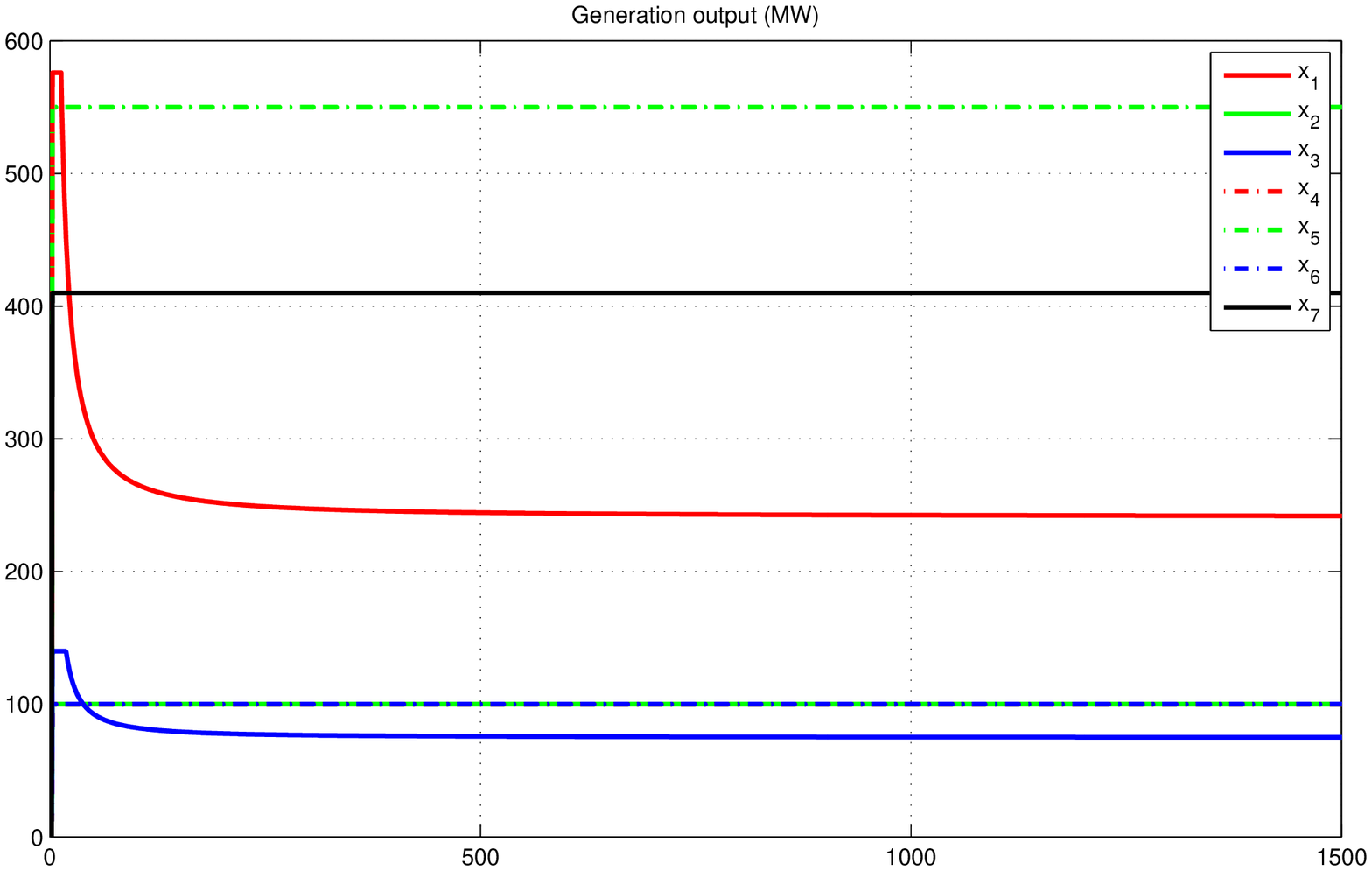}\label{fig1:subfig:b}}
\end{figure}

\begin{figure}[!htbp]
\centering
\subfigure[Total generation v.s. demand (MW))]{\includegraphics[width=18cm]{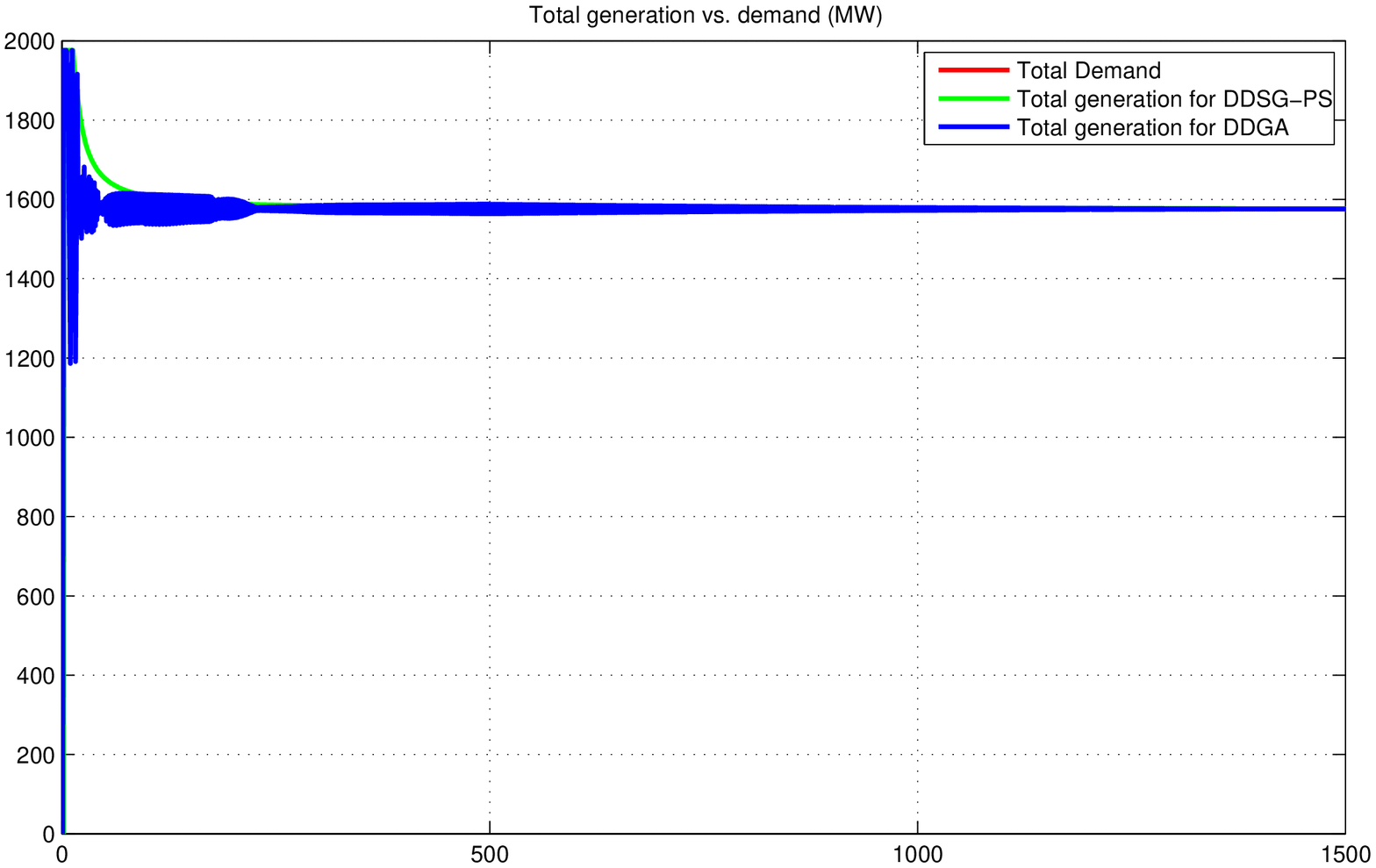}\label{fig1:subfig:c}}
\subfigure[Values of total cost function v.s. iterations ]{\includegraphics[width=18cm]{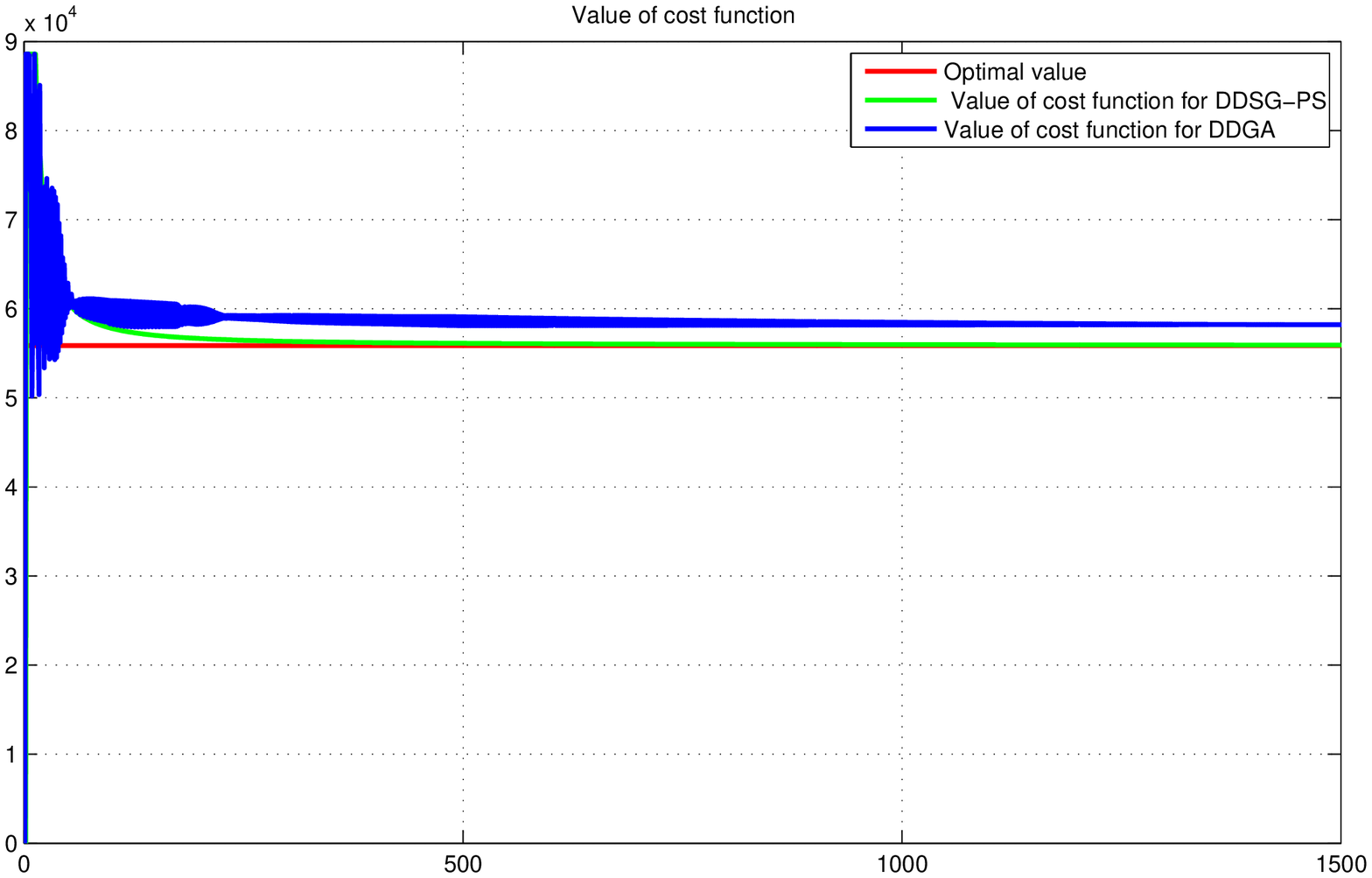}\label{fig1:subfig:d}}
\end{figure}

Figure \ref{fig1:subfig:a} shows the evolution of the values of dual variable (corresponding the incremental cost or price) at the first 30 iterations and at the first 1500 iterations. We can observe that all local dual variables $\lambda_{i}, i=1,2,\ldots, 7$ agree on the same value after earlier 50 iterations. Figure  \ref{fig1:subfig:b} demonstrates the trendy of power generations at the first 1500 iterations.  It can seen that Algorithm DDSG-PS can gradually approximate the optimal solution in a very short time.
Figure \ref{fig1:subfig:c} illustrates the evolution of total generations versus total demand. From Figure \ref{fig1:subfig:c}, we can find that the outputs of total generation indeed meet the total demand. Figures \ref{fig1:subfig:a}, \ref{fig1:subfig:b} and \ref{fig1:subfig:c} validate our theoretical results.
As shown in Figure \ref{fig1:subfig:d}, the iterative values of of total cost function are rapidly convergent to the optimal value (red solid line).

\section{Conclusions}\label{sec6}

In this paper, a distributed algorithm for convex optimization with local and coupling constraints over time-varying directed networks was proposed. The algorithm incorporated the push-sum protocol into dual subgradient methods. The optimality of primal and dual variables, and constraint violations was established. Moreover,  the explicit convergence rates of the proposed algorithm were obtained. Some numerical results showed that the proposed method is efficacy.

\end{document}